\documentclass[a4paper,12pt]{article}
\usepackage[english]{babel} 
\usepackage{amsfonts,amssymb,amsmath,a4wide,amsbsy} 
\usepackage{amsthm} 
\renewcommand{\P}{{\mathbb{P}}}

\newcommand{\thefont}[2]{\fontsize{#1}{#2}\fontshape{n}\selectfont}
\def\ind{\rlap{\thefont{10pt}{12pt}1}\kern.16em\rlap{\thefont{11pt}{13.2pt}1}\kern.4em}

\providecommand{\be}{\begin{equation}}
\providecommand{\ee}{\end{equation}}

\theoremstyle{plain}
\newtheorem{thm}{Theorem}[section]
\newtheorem{lem}[thm]{Lemma}

\theoremstyle{definition}

\newtheorem{rem}[thm]{Remark}

\title{On a conjecture concerning the sum of independent Rademacher random 
variables\footnote{The problem studied in this paper has been communicated to the
author  by Vidmantas Bentkus in the beginning of the year 2010 and this paper is 
devoted to his memory.}}
\author{Martien C.A. van Zuijlen}
\date{}

\begin{document}

\maketitle
\centerline{IMAPP, MATHEMATICS} 
\medskip
\centerline{RADBOUD UNIVERSITY  NIJMEGEN}
\begin{center}
	{
	Heyendaalseweg 135\\ 
	6525 AJ Nijmegen\\ 
	The Netherlands\\
	e-mail author: M.vanZuijlen@science.ru.nl\\
	
	
	}
\end{center}

\begin{abstract}
\begin{scriptsize}
\begin{large}
•\begin{Large}
\begin{scriptsize}
•\begin{large}\begin{small}
\begin{Huge}
•\begin{tiny}
•
\end{tiny}
\end{Huge}
\end{small}
•
\end{large}
\end{scriptsize}
\end{Large}
\end{large}
\end{scriptsize}

  For a fixed unit vector $a=(a_1,a_2,...,a_n)\in S^{n-1},$ we consider the $2^n$ sign 
 vectors $\epsilon=(\epsilon_1,\epsilon_2,...,\epsilon_n)\in \{-1,1\}^n$ and the 
 corresponding scalar products $a.\epsilon=\sum_{i=1}^n a_i\epsilon_i.$ 
 In \cite{HK} the following old conjecture has been reformulated, stating 
 that of the $2^n$ sums of the form $\sum\pm a_i$ it is impossible that there are
  more with $|\sum_{i=1}^n \pm a_i|>1$ than there  are with 
  $|\sum_{i=1}^n \pm a_i|\leq1.$ The result is of interest in itself, but has also 
   an appealing reformulation in probability theory and in geometry.
 In this paper we will solve this problem in the uniform case where all the $a$'$s$ are equal. More precisely, for $S_n$ being a sum of $n$ independent Rademacher random variables,  we will show that 
$$P_n:=\P\{-\sqrt{n}\leq S_n \leq \sqrt{n}\}\geq 1/2.$$ Hence, there is already $50\%$ of the probability mass between minus one and one standard deviation. This lower bound is sharp, it is much better than for instance the bound that can be obtained  from application of the Chebishev inequality and our bound will have nice applications in random walk theory. 

The method of proof is of interest in itself and useful for obtaining exact
 expressions and very precise estimations for the tail probabilities $P_n$ for
  finite  samples. For instance, it turns out that already for $n>7$ the tail probabilities are above $4719/8192\approx 0.58,$  whereas according to the central limit theorem the limit for these probabilities is $\approx 0.68.$
\end{abstract} 

\section{Introduction and result}
	For a positive integer $n,$ let $\epsilon_1,\epsilon_2,...,\epsilon_n$ be i.i.d. Rademacher 
random variables and let $a=(a_1,a_2,...,a_n)$ be a vector of positive 
real numbers with 
$\sum_{i=1}^na_i^2=1$. The following unsolved problem has been presented in \cite{G} 
and is attributed to B. Tomaszewski. In \cite{HK}, Conjecture 1.1, this old 
conjecture has been reformulated:
  $$\P(|a_1\epsilon_1+a_2\epsilon_2+...+a_n\epsilon_n|\leq 1)\geq 1/2.$$ 
  This problem is at least  25 years old and even in the uniform case where,
    $$a_1=a_2=...=a_n=n^{-1/2},$$   to our best knowledge the problem seems still
     to be unsolved. 
In the uniform case  the maximum possible value of $\frac{S_n}{\sqrt n}$ is 
$\sqrt n,$ where 
\begin{equation}\label{vgl0}
 S_n:=\epsilon_1+\epsilon_2+...+\epsilon_n
\end{equation}
 and the conjecture states that
\begin{equation}\label{vgl1}
 P_n:=\P\{|S_n|\leq \sqrt{n}\}=\P\{|\sum_{i=1}^n\epsilon_i|\leq \sqrt n\}\geq 1/2. 
 \end{equation}
 The present paper solves the conjecture in this special case. It means that  at least $50\%$ of the probability mass is between minus one and one standard deviation from the mean, which is quite remarkable. We note that 
 \begin{itemize}
 \item[i)] the problem can be easily reformulated in terms of Binomial(n,1/2) random variables since $\eta_i:=(\epsilon_i+1)/2, i=1,2,...,n,$ are independent Bernoulli(1/2) random variables and hence  
 \begin{equation}\label{vgl111}
 P_n:=\P\{|S_n|\leq \sqrt{n}\}=\P\{|\sum_{i=1}^n\epsilon_i|\leq \sqrt n\}=
  \P\{ \frac{n-\sqrt{n}}{2}\leq  T_n\leq \frac{n+\sqrt{n}}{2}\},
 \end{equation}
 where $T_n=\eta_1+\eta_2+...+\eta_n$ is a Binomial(n,1/2) random variable. Note that $2T_n=S_n+n.$
 \item[ii)] easy calculations show that the  sequence  $(P_n)$ is 
not monotone in $n$;  
\item[iii)] the method we use is also suitable to obtain sharp lower bounds and upper bounds for probabilities concerning $a$ standard deviations ($a>0)$:
$$\P\{|S_n|\leq a\sqrt n\}= 
\P\{|\sum_{i=1}^n\epsilon_i|\leq a\sqrt n\}.$$
\end{itemize}

 The result is as follows.

\begin{thm}
 Let $ \epsilon_1,\epsilon_2,...,\epsilon_n$ be independent Rademacher random variables, so that $$\P\{\epsilon_1=1\}=\P\{\epsilon_1=-1\}=1/2$$ and let
  $S_n$  and 
 $P_n$ be defined as in (\ref{vgl0}) and (\ref{vgl1}).
 Moreover, let $$N:=\{0,1,2,3,...\}=\bigcup_{k=1}^{\infty} C_k,$$ 
where 
$$C_k=\{n: k^2\leq n+1<(k+1)^2\}=$$
\begin{equation}\label{vgl2}
=\{k^2-1,k^2,k^2+1,...,(k+1)^2-3,(k+1)^2-2\}.
\end{equation}
Then,  for $k=1,2,...,$ we have
\begin{equation}\label{vgl3}
   Q_k^-:=P_{(k+1)^2-2}=\min_{n\in C_k}P_n,
   \end{equation}
  \begin{equation}\label{vgl31}
   Q_k^+:=P_{k^2}=\max_{n\in C_k}P_n.
   \end{equation} 
  Moreover, the sequence $(Q_k^-)$ is monotone increasing, the sequence $(Q_k^+)$ is monotone decreasing and
  \begin{equation}\label{vgl32}
  lim_{k\rightarrow \infty} Q_k^- =lim_{k\rightarrow \infty} Q_k^+ = \Phi(1)\approx 0.68, 
  \end{equation}
  where $\Phi$ is the standard normal distribution function.
  It follows that

 $$ P_n \begin{cases} =1, & \;\;\text{for} \;\; n=1,\\
 			=\frac{1}{2}, & \;\;\text{for} \;\; n=2,\\
  			 \geq \frac{35}{64}\approx 0.55, & \;\;\text{for}\;\; n=3,4,5,6,7,\\
	          \geq \frac{4719}{8192}\approx0.58, & \;\;\text{for}\;\; n=8,9,...,14,\\
	            \geq \frac{156009}{262144}\approx 0.60, & \;\;\text{for}\;\; n\geq 15.\\
	            \end{cases} $$ 
 so that certainly  for $n\in N$
 $$P_n\geq 1/2,$$
  and 
  $$ P_n \begin{cases} =1, & \;\;\text{for} \;\; n=1,\\
 			=\frac{1}{2}, & \;\;\text{for} \;\; n=2,\\
  			 \leq \frac{7}{8}\approx 0.88, & \;\;\text{for}\;\; n=3,4,5,6,7,\\
	          \leq \frac{105}{128}\approx 0.82, & \;\;\text{for}\;\; n=8,9,...,14,\\
	            \leq \frac{25833}{32768}\approx 0.79, & \;\;\text{for}\;\; n\geq 15.\\
	             \leq \frac{3231615}{4194304}\approx 0.77, & \;\;\text{for}\;\; n\geq 24.\\
	            \end{cases}. $$ 
 so that certainly  for $n\geq 2$
 $$P_n\leq 7/8.$$
 Note that trivially $P_1=1,$  $P_2=1/2$  and by definition $P_0=1.$
 \end{thm}
 
 \begin{rem} It will turn out that our method of proof  also gives a
 way to compute $P_n,$ for $n=1,2,...,$ in a recursive way.
 More precisely,
 let $n\in C_{k_n}$ for some $k_n\geq 2$ and $n=k_n^2-1+i$ for some $i\in \{0,1,...,2k_n\}$, 
 then

$$P_n=1/2 +$$
$$+\sum_{k=2}^{k_n-1}\left\{\P\{S_{k^2-2}=k\}
+\sum_{j=1,j=odd}^{2k}\P\{S_{k^2+j-2}=k+1\}-
\sum_{j=1,j=even}^{2k}\P\{S_{k^2+j-2}=k\}\right\}+$$
$$ +\P\{S_{k_n^2-2}=k_n\}+\sum_{j=1,j=odd}^{i}\P\{S_{k^2+j-2}=k_n+1\}-
\sum_{j=1,j=even}^{i}\P\{S_{k^2+j-2}=k_n\}.$$

 \end{rem}
 \section{Proof of the Theorem}

 Note that because of the symmetry in the distribution of $S_n$ we have
 $$P_n:=\P\{-\sqrt{n}\leq S_n \leq \sqrt{n} \}=\P\{1-\sqrt{n}\leq S_{n-1} \leq 1+\sqrt{n} \}.$$
 We will  compare the probability 
 $P_n$ with
 $$P_{n-1}:=\P\{-\sqrt{n-1}\leq S_{n-1} \leq \sqrt{n-1} \}.$$
 Notice that $S_n$ takes on values with positive probability only in the 
 set $D_n,$ where $D_n=\{0,2,-2,4,-4,...,n,-n\}$ for $n$ is even and $D_n=\{1,-1,3,-3,...,n,-n\}$ 
 for $n$ is odd. Also note that the different values of $D_n$ have at least a distance $2$, so that there can be at most one element of $D_n$ in a half open interval of lenght $2.$ 
  In fact $$P_n=P_{n-1}-\P\{S_{n-1}\in B_{n}\}+\P\{S_{n-1}\in \tilde{A}_{n}\}$$
  \begin{equation}\label{vgl4}
  =P_{n-1}-\P\{S_{n-1}\in B_{n}\}+\P\{S_{n-1}\in A_{n}\},
  \end{equation}
  where
 $$B_{n}=[-\sqrt{n-1},1-\sqrt{n}),$$  $$\tilde{A}_{n}=(\sqrt{n-1},1+\sqrt{n}]$$ and $$A_{n}=[-1-\sqrt{n},-\sqrt{n-1})
 .$$ Since   
   $$ 0<\sqrt{n}-\sqrt{n-1}< 1$$   the difference 
 between $P_n$ and $P_{n-1}$ will be rather small and  only some integers have to be taken into account. In fact, it follows that
 $$|B_n|=|1-(\sqrt{n}-\sqrt{n-1})|\in (0,1), \;\;  \;\; |A_n|=|1+(\sqrt{n}-\sqrt{n-1})|\in (1,2)\;\;\; \;\; $$
 and
 $$|A_n\cup B_n|=|[-1-\sqrt{n},1-\sqrt{n})|=2,$$
 where $|C|$ is the length of the interval $C.$ It is clear that
  $A_n\cup B_n$  contains exactly one even integer and  one odd integer
and the set $B_{n}$ contains at most one
 negative integer and the set $A_{n}$ contains at most two negative  integers. If we restrict ourselves to the even integers, then 
  exactly one of the sets $A_n$ and $B_n$ contains 
  one even integer. The other set contains no even integer.
The same statement holds if we restrict ourselves to the odd integers. 
 
The main problem is that we have to show that $P_n\geq 1/2,$ for $n\in \{2,3,4,...\}.$ We will solve the problem
separately for blocks of integers related with two adjacent squares. Note that
$$N=\{0,1,2,3,...\}=\bigcup_{k=1}^{\infty} C_k,$$ 
where 
\begin{equation}\label{vgl5}
C_k=C_{k,1}\cup C_{k,2},
\end{equation}
and
$$C_{k,1}:=\{k^2,k^2+2,...,(k+1)^2-3\}=
\begin{cases}\{n\in C_k:\;\; n\;\; \text{is even}\},\;\;\text{if}\;\;k \;\;\text{is even},\\
\{n\in C_k: \;\;n\;\; \text{is odd}\},\;\;\text{if}\;\;k \;\;\text{is odd},\\
\end{cases}$$
$$C_{k,2}:=\{k^2-1,k^2+1,...,(k+1)^2-2\}=
\begin{cases}\{n\in C_k:\;\; n\;\; \text{is even}\},\;\;\text{if}\;\;k \;\;\text{is odd},\\
\{n\in C_k: \;\;n\;\; \text{is odd}\},\;\;\text{if}\;\;k \;\;\text{is even}.\\
\end{cases}$$
Note that $C_k$ contains $2k+1$ non-negative integers, $C_{k,1}$ contains $k$ non-negative integers, and $C_{k,2}$ contains $k+1$ non-negative integers. In case $k$ is even (odd),  then $C_{k,1}$ consists of all even (odd) integers in $C_k$ and $C_{k,2}$ consists of all odd (even) integers in $C_k,$ since 
$$[k=\text{even}]\Leftrightarrow [k^2-1=\text{odd}].$$
 Furthermore, we  have
  $$\{-k,1-k\}\in  A_n\cup B_{n}= [-1-\sqrt n,1-\sqrt n),\;\; \text{for}\;\;n=k^2-1$$
 and
 $$\{-1-k,-k\}\in  A_n\cup B_{n}= [-1-\sqrt n,1-\sqrt n),\;\; \text{for}\;\;n\in C_k\backslash\{k^2-1\}.$$
Moreover,  for $n\in C_k$ we have
 $$\P\{S_{n-1}\in A_n\cup B_{n}\}=\P\{S_{n-1}\in [-1-\sqrt n,1-\sqrt n)\}=
\begin{cases}\P\{S_{n-1}=-1-k\},\;\; \text{for} \;\; n\in C_{k,1},\\   
\P\{S_{n-1}=-k\}, \;\; \text{for} \;\; n\in C_{k,2}.
\end{cases}$$  
For each  $n\in C_k$ by using (\ref{vgl4}) we can compute $P_n$ from $P_{n-1}$ by subtracting $\P\{S_{n-1}\in B_{n}\}$ or adding
$\P\{S_{n-1}\in A_{n}\}.$ It will turn out that for the integers $n\in C_{k,2}\backslash
\{k^2-1\}$   we have negative contributions due to $\P\{S_{n-1}\in B_n\},$ whereas for the
integers $n\in C_{k,1}\cup \{k^2-1\} $ we have positive contributions due to $\P\{S_{n-1}\in A_n\}.$ For instance, for
$n=3$ we obtain
$$\P\{S_{2}\in \tilde{A}_3\}=\P\{S_{2}\in A_3\}=\P\{\epsilon_1+\epsilon_2\in [-1-\sqrt{3},-\sqrt 2)\}=\P\{\epsilon_1+\epsilon_2=-2\}={1/2}^2=1/4,$$
 for $n=4$ we obtain
$$\P\{S_{3}\in \tilde{A}_4\}=\P\{S_{3}\in A_4\}=\P\{\epsilon_1+\epsilon_2+\epsilon_3\in [-3,-\sqrt 3)\}=\P\{\epsilon_1+\epsilon_2+\epsilon_3=-3\}={1/2}^3=1/8,$$
and for $n=7$ we obtain
$$\P\{S_{6}\in B_7\}=\P\{\epsilon_1+\epsilon_2+...+\epsilon_6\in [-\sqrt 6,1-\sqrt7)\}=\P\{\epsilon_1+\epsilon_2+...+\epsilon_6=-2\}=\binom{6}{2}2^{-6}.$$ More generally, we obtain the following lemma:

\begin{lem}\label{lem1} For $k=2,3,...$ we have the following.

\begin{itemize}
\item[i)] For $n=k^2-1$ we have  $$-k\in A_n=[-1-\sqrt{k^2-1},-\sqrt{k^2-2}),$$ 
so that 

$$\P\{S_{n-1}\in A_n\}= \P\{S_{k^2-2}=k\}=\frac{\binom{k^2-2}{k(k-1)/2-1}}{2^{k^2-2}};\;\;\;
\P\{S_{n-1}\in B_n\}=0.$$

$$ $$
 \item[ii)]For $n\in C_{k,1}$ we have for $n=k^2+2i,$ with $i=0,1,...,k-1,$
$$-1-k\in A_n=[-1-\sqrt{k^2+2i},-\sqrt{k^2+2i-1}),$$   so that 

$$\P\{S_{n-1}\in A_n\}= \P\{S_{k^2+2i-1}=k+1\}=\frac{\binom{k^2+2i-1}{k(k-1)/2+i-1}}{2^{k^2+2i-1}};\;\;\;
\P\{S_{n-1}\in B_n\}=0.$$

$$ $$
\item[iii)]
For $n\in C_{k,2}\backslash\{k^2-1\}$ we have for $n=k^2+1+2i,$ with $i=0,1,...,k-1,$
 $$-k\in B_n=[-\sqrt{k^2+2i},1-\sqrt{k^2+1+2i}),$$   so that 
 
$$\P\{S_{n-1}\in B_n\}= \P\{S_{k^2+2i}=k\}=\frac{\binom{k^2+2i}{k(k-1)/2+i}}{2^{k^2+2i}};\;\;\;
\P\{S_{n-1}\in A_n\}=0.$$
\end{itemize}

\end{lem}
$$ $$
 For $k=2,3,4$ and also for general $k\in \{2,3,...,\}$ the results of the lemma above are illustrated in  the  Tables 2-5.
  The starting points where $n=1$ and $n=2$ are indicated in Table 1 separately. The proof of the lemma  will be given in Section 3.
 $$ $$

 $$ $$ 

\begin{table}
\begin{center}
 \begin{tabular}{|c|c|c|c|}
\hline $n$ & $\P\{S_{n-1}\in A_n\}$ & $\P\{S_{n-1}\in B_n\}$ & $P_n$   \\
 \hline $0$ &  &   & $-$ \\ 
\hline $1$ &  &   & $1$ \\ 
\hline $2$ &  &   & $0.50$ \\ 
\hline 
\end{tabular}
\caption{Table for $C_1$}
\end{center}
\end{table}

$$ $$
\begin{table}
\begin{center}
 \begin{tabular}{|c|c|c|c|c|}
\hline $n$ & $\P\{S_{n-1}\in A_n\}$ & $\P\{S_{n-1}\in B_n\}$ & $P_n$    \\ 
\hline $3$ &  $\frac{\binom{2}{0}}{2^{2}}$ & & $\frac{3}{4}=0.75 $\\ 
\hline $4$ & $\frac{\binom{3}{0}}{2^{3}}$ &  & $\frac{7}{8}\approx 0.88 $\\ 
\hline $5$ &   &$\frac{\binom{4}{1}}{2^{4}}$ & $\frac{5}{8}\approx 0.63 $\\ 
\hline $6$ & $\frac{\binom{5}{1}}{2^{5}}$ &  & $\frac{25}{32} \approx 0.78$ \\ 
\hline $7$ &    &$\frac{\binom{6}{2}}{2^{6}}$ & $\frac{35}{64} \approx 0.55$ \\ 
\hline 
\end{tabular}
\caption{Table for $C_2$}
\end{center}
\end{table}

$$ $$

\begin{table}
\begin{center}
 \begin{tabular}{|c|c|c|c|}
\hline $n$ & $\P\{S_{n-1}\in A_n\}$ & $\P\{S_{n-1}\in B_n\}$ & $P_n$ \\ 
\hline $8$ & $\frac{\binom{7}{2}}{2^{7}}$ & & $\frac{91}{128}\approx 0.71 $ \\ 
\hline $9$ & $\frac{\binom{8}{2}}{2^{8}}$ & & $\frac{105}{128}\approx 0.82 $ \\ 
\hline $10$ &  & $\frac{\binom{9}{3}}{2^{9}}$ & $\frac{21}{32}\approx 0.66$ \\ 
\hline $11$ & $\frac{\binom{10}{3}}{2^{10}}$ &  & $\approx 0.77$ \\ 
\hline $12$ &    &$\frac{\binom{11}{4}}{2^{11}}$ & $\approx 0.61$\\ 
\hline $13$ & $\frac{\binom{12}{4}}{2^{12}}$ &     & $\approx 0.73$  \\ 
\hline $14$ &    &$\frac{\binom{13}{5}}{2^{13}}$ & $\approx 0.58$ \\
\hline 
\end{tabular}
\caption{Table for $C_3$}
\end{center}
\end{table}

$$ $$

\begin{table}
\begin{center}
\begin{tabular}{|c|c|c|c|}
\hline $n$ & $\P\{S_{n-1}\in A_n\}$ & $\P\{S_{n-1}\in B_n\}$ & $P_n$ \\ 
\hline $15$ & $\frac{\binom{14}{5}}{2^{14}}$ &  &$\approx 0.70$ \\ 
\hline $16$ &$\frac{\binom{15}{5}}{2^{15}}$  & & $\approx 0.79$\\ 
\hline $17$ &   & $\frac{\binom{16}{6}}{2^{16}}$& $\approx 0.67$\\ 
\hline $18$ & $\frac{\binom{17}{6}}{2^{17}}$ &  & $\approx 0.76$\\ 
\hline $19$ &    &$\frac{\binom{18}{7}}{2^{18}}$& $\approx 0.64$\\ 
\hline $20$ & $\frac{\binom{19}{7}}{2^{19}}$ &  & $\approx 0.74$\\ 
\hline $21$ &    &$\frac{\binom{20}{8}}{2^{20}}$& $\approx 0.62$\\
\hline $22$ &   $\frac{\binom{21}{8}}{2^{21}}$& & $\approx 0.71$\\ 
\hline $23$ & &  $\frac{\binom{22}{9}}{2^{22}}$  & $\approx 0.60$\\
\hline 
\end{tabular}
\caption{Table for $C_4$}
\end{center}
\end{table}

$$ $$
\begin{table}
\begin{center}
 \begin{tabular}{|c|c|c|}
\hline $n$ & $A_n$ & $B_n$   \\ 
\hline $k^2-1$ & $[-1-\sqrt{k^2-1},-\sqrt{k^2-2})$&   \\ 
\hline $k^2$ & $[-1-k,-\sqrt{k^2-1})$ & \\ 
\hline $k^2+1$ & &  $[-k,1-\sqrt{k^2+1})$  \\ 
\hline $k^2+2$ & $[-1-\sqrt{k^2+2},-\sqrt{k^2+1})$&    \\ 
\hline $k^2+3$ &  & $[-\sqrt{k^2+2},1-\sqrt{k^2+3})$    \\ 
\hline $k^2+4$  & $[-1-\sqrt{k^2+4},-\sqrt{k^2+3})$ &  \\ 
\hline $...$ & $...$ & $...$    \\ 
\hline $...$ & $...$ & $...$    \\ 
\hline $...$ &$...$ & $...$     \\
\hline $...$ & $...$ & $...$    \\ 
\hline $(k+1)^2-3$ &  $[-1-\sqrt{(k+1)^2-3},-\sqrt{(k+1)^2-4})$ &   \\ 
\hline $(k+1)^2-2$  & & $[-\sqrt{(k+1)^2}-3,1-\sqrt{(k+1)^2-2})$  \\

\hline 
\end{tabular}
\caption{Intervals for $C_k$}
\end{center}
\end{table}

$$ $$
\begin{table}
\begin{center}
 \begin{tabular}{|c|c|c|c|}
\hline $n$  & $integer$ & $\P\{S_{n-1}\in A_n\}$ & $\P\{S_{n-1}\in B_n\}$  \\ 
\hline $k^2-1$ & $-k\in A_n$  &   $\P\{S_{k^2-2}=k\}=\frac{\binom{k^2-2}{k(k-1)/2-1}}{2^{k^2-2}}$& \\ 
\hline $k^2$ & $-1-k\in A_n$  &   $\P\{S_{k^2-1}=k+1\}=\frac{\binom{k^2-1}{k(k-1)/2-1}}{2^{k^2-1}}$& \\ 
\hline $k^2+1$ &$-k\in B_n$ &  &$\P\{S_{k^2}=k\}=\frac{\binom{k^2}{k(k-1)/2}}{2^{k^2}}$  \\ 
\hline $k^2+2$ & $-1-k\in A_n$ &   $\P\{S_{k^2+1}=k+1\}=\frac{\binom{k^2+1}{k(k-1)/2}}{2^{k^2+1}}$& \\ 
\hline $k^2+3$ & $-k\in B_n$  & &  $\P\{S_{k^2+2}=k\}=\frac{\binom{k^2+2}{k(k-1)/2+1}}{2^{k^2+2}}$ \\ 
\hline $k^2+4$ &$-1-k\in A_n$  &   $\P\{S_{k^2+3}=k+1\}=\frac{\binom{k^2+3}{k(k-1)/2+1}}{2^{k^2+3}}$ & \\ 
\hline $...$ & $...$ &   $...$ &\\ 
\hline $...$ &  $...$ &  & $...$  \\ 
\hline $...$ & $...$ &     &\\
\hline $...$ &  $...$ &   $...$ &\\ 
\hline $(k+1)^2-3$ &$-1-k\in A_n$     & $\P\{S_{(k+1)^2-4}=k+1\}=\frac{\binom{(k+1)^2-4}{k(k+1)/2-2}} {2^{(k+1)^2-4}}$& \\ 
\hline $(k+1)^2-2$ &$-k\in B_n$  &&
$\P\{S_{(k+1)^2-3}=k\}=\frac{\binom{(k+1)^2-3}{k(k+1)/2-1}} {2^{(k+1)^2-3}}$  \\

\hline 
\end{tabular}
\caption{Probabilities for $C_k$}
\end{center}
\end{table}

$$ $$
$$ $$ 



 

 
 To complete the proof one has to show that the other $P_n$ for $n\in C_k$ cannot be smaller
than $P_{(k+1)^2-2}.$ This can be done by comparing a contribution due to $B_n$ with the
lower but adjacent contribution due to $A_n:$ 
 \begin{lem}\label{lem2}
 For  $i=0,1,...,k-1,$ let
 \begin{equation}\label{vgl11}
  \delta_i:=\P\{S_{k^2+2i-1}=k+1\}- \P\{S_{k^2+2i}=k\}=
\frac{\binom{k^2-1+2i}{k(k-1)/2+i-1}}{2^{k^2-1+2i}}-\frac{\binom{k^2+2i}{k(k-1)/2+i}}{2^{k^2+2i}}.
  \end{equation} 
  Then we have
   \begin{equation}\label{vg22}
 \delta_0\leq \delta_1\leq ...\leq\delta_{k-1}<0. 
 \end{equation}
 \end{lem}
The proof of this lemma  will be given in Section 3.
 $$ $$
 The remaining line of proof of Theorem 1.1 is now as follows.
   Because of the monotony of the $\delta$'s we have
    $$ [\P\{S_{k^2-2}=k\}+k\delta_0\geq 0]\Rightarrow [\P\{S_{k^2-2}=k\}+\sum_{i=0}^{k-1}\delta_i\geq 0]\Leftrightarrow [P_{k^2-2}\leq P_{(k+1)^2-2}]$$
    and
    $$ \P\{S_{k^2-2}=k\}+k\delta_0$$
    $$=\frac{\binom{k^2-2}{k(k-1)/2-1}}{2^{k^2-2}}-k\binom{k^2-1}{k(k-1)/2-1}\frac{2k}{2^{k^2}(k(k-1))}=0,$$
    so that we obtain indeed
    $$P_{k^2-2}\leq P_{(k+1)^2-2}. $$
  Next, using that the $\delta$'s are negative, we have   for $k=2,3,...,$ 
$$Q_k^-:=P_{(k+1)^2-2}=\min_{n\in C_k}P_n, $$ since 
$$P_{(k+1)^2-2}\leq P_{(k+1)^2-4}\leq ...\leq P_{k^2-1} $$ and also
$$P_{(k+1)^2-2}\leq P_{(k+1)^2-3}\leq P_{(k+1)^2-5}\leq ...\leq P_{k^2}.$$ 
Since trivially $ P_{k^2-1}\leq P_{k^2},$ it also follows that
$$Q_k^+:=P_{k^2}=\max_{n\in C_k}P_n, $$ 
which completes the proof of the theorem.
 
  \section{Proof of the Lemmas}

  In this section we will deliver the proof of Lemma  \ref{lem1}  and the proof of Lemma  \ref{lem2}.
\\  
\\
\textbf{Proof of Lemma \ref{lem1}}:
  Suppose for the moment that  $n-1$ is even so that the values of $S_{n-1}$ can be only even
  integers with positive probability. Since $|A_n\cup B_n|=2$ and since the interal $A_n\cup
  B_n$ is half open, there will be exactly one even integer in $A_n\cup B_n,$ so that either 
  $\P\{S_{n-1}\in A_n\}=0$ or  $\P\{S_{n-1}\in B_n\}=0.$ The same argument applies in case $n-1$ is odd.
  
  To prove the first part of the lemma, let   $n=k^2-1.$ Since
  $$-1-k<-1-\sqrt{k^2-1}\leq -k <-\sqrt{k^2-2}$$ we have
  \begin{equation}\label{vgl12}
  -k\in A_n=[-1-\sqrt{k^2-1},-\sqrt{k^2-2}).
  \end{equation}
  Moreover, there is exactly one possible value of $S_{n-1}$ in the interval $$A_n\cup
  B_n=[-1-\sqrt{k^2-1},1-\sqrt{k^2-1})$$ and since $n-1=k^2-2$ and $k$ are either both even or
  both odd, the only possible value of $S_{n-1}$ in the interval $A_n\cup B_n$ is $-k\in A_n,$ which implies that 
  $\P\{S_{n-1}\in B_n\}=0$
  and 
   $\P\{S_{n-1}\in A_n\}=\P\{S_{n-1}=k\}.$
   This last probability can easily be calculated by using (\ref{vgl111}).
   
   Next, let $n\in C_{k,1},$ so that $n=k^2+2i,$ for some $i\in \{0,1,...,k-1\}.$ Note that 
    $$-1-\sqrt{k^2+2i}\leq -1-k\leq -k <1-\sqrt{k^2+2i},$$ so that 
    $$\{-1-k,-k\}\subset A_n\cup B_n=[-1-\sqrt{k^2+2i},1-\sqrt{k^2+2i}). $$
   Moreover,  since $n-1=k^2+2i-1,$ and $1+k$ are either both even or both odd, the only
   possible value of $S_{n-1}$ in the interval $A_n\cup B_n$ is $-1-k\in A_n,$ which implies that 
  $\P\{S_{n-1}\in B_n\}=0$
  and 
   $\P\{S_{n-1}\in A_n\}=\P\{S_{n-1}=1+k\}.$
   This last probability can easily be calculated again.

   Finally, for $n\in C_{k,2}\backslash \{k^2-1\},$ a similar reasoning can be followed.

 $$ $$ 
 
 \noindent \textbf{Proof of Lemma \ref{lem2}}:
We will use for $i=0,1,...,k-1$ the identity
 $$\binom n j 2^{-n}-\binom{n+1}{j+1}2^{-n-1}=\binom n j 2^{-n-1}\{\frac{2j-n+1}{j+1}\}, $$
 with $n=k^2-1+2i$ and $j=k(k-1)/2+i-1,$ and obtain 
 $$\delta_i=-\frac{\binom{k^2-1+2i}{k(k-1)/2+i-1}}{2^{k^2+2i}}\frac{2k}{k(k-1)+2i}<0. $$
 Moreover, for $i\in\{0,1,...,k-2\},$  we find 
 $$[\delta_i \leq \delta_{i+1}]\Leftrightarrow $$
 
$$\Leftrightarrow  [\frac{\binom{k^2-1+2i}{k(k-1)/2+i-1}}{2^{k^2+2i}}\frac{2k}{k(k-1)+2i}\geq \frac{\binom{k^2-1+2i+2}{k(k-1)/2+i}}{2^{k^2+2i+2}}\frac{2k}{k(k-1)+2i+2}]$$
 
 $$\Leftrightarrow [\frac{4}{k(k-1)+2i}\geq \frac{(k^2+2i)(k^2+2i+1)}{\{k(k-1)/2+i\}\{k(k+1)/2+i+1\}\{(k(k-1)+2i+2\}}]$$
 
 $$\Leftrightarrow [2+3i+k^2\geq 0].$$
 This last inequality holds trivially, so that (\ref{vg22}) holds.

 $$ $$
 
\noindent \textbf{Acknowledgement}. The author is indebted to Harrie Hendriks for careful 
reading the manuscript and for giving useful comments, which have led to a
 simplification in the proof of the main theorem.

\end{document}